\documentclass[11pt]{article}
\usepackage{amsfonts,amsmath,enumerate}

\newtheorem{theorem}{Theorem}

\def\C{{\mathbb C}}
\def\R{{\mathbb R}}
\def\D{{\cal D}}
\def\W{{\cal W}}
\def\Z{{\mathbb Z}}

\begin{document}

\title{Spectral curves for Cauchy--Riemann operators on punctured
  elliptic curves} \author{C.\ Bohle \thanks{Universit\"at T\"ubingen,
    Fachbereich Mathematik, Auf der Morgenstelle 10, 72076 T\"ubingen,
    Germany; e-mail:bohle@mathematik.uni-tuebingen.de} \and I.A.\
  Taimanov \thanks{Sobolev Institute of Mathematics, 630090
    Novosibirsk, Russia; e-mail: taimanov@math.nsc.ru.} } \date{}

\maketitle

{\bf 1.\ Introduction}. \ \ The spectral curve of a periodic
differential operator is fundamental to finite gap integration
theory. The solution of the inverse problem for a finite gap operator
consists in assigning such an operator to the spectral data consisting
of the spectral curve and some additional data related to it. For
one--dimensional Schr\"odinger operators the spectral curve was defined
in \cite{Novikov1974}, and in \cite{DKN} this notion was generalized to
two--dimensional Schr\"odinger operators for which the spectral curve
parametrizes the
Floquet--Bloch functions at a fixed energy level.

In the present article we show that an analogous spectral curve can be
defined for the Cauchy--Riemann operator on a punctured elliptic
curve, provided one imposes appropriate boundary conditions.

\begin{theorem}
  Let $X$ be an elliptic curve $X = \C/\{\Z e_1 + \Z e_2\}$, where
  $e_1$ and $e_2$ are generators of the period lattice, and let $p_1$,
  \dots, $p_N$ be pair--wise distinct points on $X$.  We consider the
  linear problem
\begin{equation}
\label{e1}
\bar{\partial}\psi = 0,
\end{equation}
\begin{equation}
\label{e2}
\psi(z+e_j)=e^{2\pi i \langle k, e_j\rangle }\psi(z), \ \ \ j=1,2,
\end{equation}
where $k = (k_1,k_2) \in \C^2, \langle k,e_j \rangle = k_1
\mathrm{Re}\,e_j + k_2 \mathrm{Im}\,e_j$ and $\psi$ satisfies
\begin{equation}
\label{e3}
\psi(z) = \frac{a_j}{z-p_j} + O(|z-p_j|) \ \mbox{for $j=1,\dots,N$}
\end{equation}
for some $a_j\in \C$.  The solutions to (\ref{e1})--(\ref{e3}) with $k$
and $a_1$, ..., $a_N$ arbitrary are parametrized (up to scale) by a
one--dimensional (spectral) curve, see (\ref{Ncurve}), which up to one puncture is
an $N$--sheeted covering of $X$.
\end{theorem}

Algebraic curves of the type thus obtained appear as irreducible components of
spectral curves of minimal tori with planar ends in $\R^3$ (see \S 4
below). It appears that
these curves coincide with the spectral curves of certain elliptic KP 
solitons as studied by Krichever \cite{Krichever1980}. 

It should be noted that, unlike in the theory of parabolic bundles, we
do not prescribe principal parts (as the $a_j$ are not fixed), but cut
down dimensions by demanding that the order zero terms in the
asymptotic expansions (\ref{e3}) at the punctures vanish.

{\bf 2.\ Spectral curves of Schr\"odinger operators.}
Let $L = \partial \bar{\partial} +U$ be a two-dimensional Schr\"odinger operator with
a double-periodic potential: $U(z+e_1)=U(z+e_2)=U(z)$.
We say that $\psi$ is a Floquet eigenfunction of $L$ with the
eigenvalue $E$ if it is a formal solution to the equation
$$
L \psi = E \psi
$$
which satisfies the periodicity conditions (\ref{e2}).  The quantities
$k_1$ and $k_2$ are called the quasimomenta of $\psi$ and
$(\nu_1,\nu_2)= (e^{2\pi i \langle k,e_1 \rangle},e^{2\pi i \langle k,e_2 \rangle})$ are the
Floquet multipliers of $\psi$.

The spectral curve which parametrizes the Floquet functions for
$E=0$ up to scale was defined in \cite{DKN}, where it was shown
how to reconstruct an operator $L$ from a (complex) curve of finite
genus and some additional data.  Such operators are called finite gap
at the zero energy level.

In general spectral curves are not finite genus and
there are two ways to rigorously establish their existence:
by using the Fredholm alternative for analytic
pencils of operators \cite{Taimanov1998}, or
by perturbing the spectral curve for
$L_0=\partial\bar{\partial}$ \cite{Krichever1989}.

% The classical definition of the spectral curve for one-dimensional Schr\"o\-din\-ger operator
% $$
% L= \frac{d^2}{dx^2} + u(x), \ \ \ u(x+T)=u(x),
% $$
% is as follows \cite{Novikov1974}. By definition, a Bloch function $\phi$ meets satisfies the conditions:
% $$
% L\phi = E\phi, \ \ \ \phi(x+T) = e^{2\pi i T}\phi(x).
% $$
% Such functions are parameterized by hyperelliptic curves, which are in general of infinite genus, of the form
% $\lambda^2 = F(E)$ which is the spectral curve of $L$. Let us consider a two-dimensional operator
% $$
% \widetilde{L} = \partial\bar{\partial}+\frac{u(x)}{4},
% $$
% then $\psi(x,y) = e^{i\sqrt{E}y}\phi(x)$ is the Floquet function of $\widetilde{L}$.

{\bf 3.\ Spectral curves of tori in $\R^3$.} \ \ The Weierstrass
representation assigns to every solution $\psi\colon
X \rightarrow \C^2$ of the problem
$$
\D \psi = 0,
$$
where $\D$ is a Dirac operator with a double--periodic real-valued potential
$$
\D =
\left(
\begin{array}{cc}
0 & \partial \\ -\bar{\partial} & 0
\end{array}
\right)
+
\left(
\begin{array}{cc}
U & 0 \\ 0 & U
\end{array}
\right), \ \ \ \ U = \bar{U},
$$
a conformal immersion of $X$ or its covering into $\R^3$,
determined by
$$
x^k = x^k(0) + \int \left( x^k_z dz + \bar{x}^k_z d\bar{z}\right), \ \ k=1,2,3,
$$
with
\begin{equation}
\label{int}
x^1_z = \frac{i}{2}(\bar{\psi}^2_2 + \psi^2_1), \ \ \
\ x^2_z = \frac{1}{2}(\bar{\psi}^2_2 - \psi^2_1), \ \ \ \ x^3_z =
\psi_1\bar{\psi}_2,
\end{equation}
branched at zeros of $\psi$ and with a double-periodic Gauss map
(see \cite{K,Taimanov1995}). Moreover, every conformal immersion of
$X$ does admit such a representation if one allows for $\psi$ to have
$\Z_2=\{\pm 1\}$--monodromy, i.e., to be the section of a spinor
bundle on $X$ \cite{Taimanov1995}.

Such a representation exists for all
immersed surfaced \cite{Taimanov1995},  and for $U=0$ it
reduces to the
classical Weierstrass representation of minimal surfaces.

The spectral curve of a torus in $\R^3$ is defined, in the spirit of
\cite{DKN}, as the spectral curve of the corresponding Dirac operator
$\D$ at the zero energy level \cite{Taimanov1998}.  It parametrizes
Floquet functions, taken up to scale, that satisfy the equation
$\D\psi=0$ and the conditions (\ref{e2}).

This spectral curve encodes the Willmore functional $\W = \int H^2
d\mu = 4 \int U^2 dx \wedge dy$, where $H$ is the mean curvature of a
closed oriented surface immersed into $\R^3$ and $d\mu$ is the measure
induced by the immersion.  Willmore tori, i.e.\ the critical points of
the Willmore functional, are described by integrable systems. They
split into several classes one of which are the  stereographic
projections of minimal tori with planar ends \cite{Bryant}.

For Willmore tori that are not Euclidean minimal with planar ends it
was proven in \cite{Bohle2010} that the spectral curve can be
compactified by adding two points at infinity.  (The result of
\cite{Bohle2010} more generally concerns constrained Willmore tori,
i.e.\ the critical points of $\W$ restricted to tori within a fixed
conformal class.)

{\bf 4.\ Spectral curves of minimal tori with planar ends.} \ \ The
spectral curves of Willmore tori coming from minimal tori with planar
ends have recently been described in \cite{BT} as reducible algebraic
curves with two irreducible components $\Gamma_+$ and $\Gamma_-$ which
are complex conjugate to each other and exactly of the type arising
in Theorem~1.

In terms of the Weierstrass representation minimal surfaces correspond
to the trivial potential $U=0$.  Hence their spinor $\psi$ splits into
two components $\psi_1,\psi_2$ such that $\psi_1$ and $\bar{\psi}_2$
satisfy
$$
\bar{\partial}\psi=0.
$$

By a planar end we mean a minimal conformal immersion $f: D \setminus
\{p\} \to \R^3$ of a punctured disc $D \setminus \{p\}$ such that
$\lim_{q \to p} f(q) = \infty$ and the immersion extends to an
immersed disc in the conformal $3$--sphere $S^3=\R \cup
\{\infty\}$.  Such end $p$ is characterized by $x^j_z, j=1,2,3$ having
poles of order $2$ with vanishing residues.  In terms of $\psi =
(\psi_1,\psi_2)$ this condition becomes (\ref{e3}).

{\bf 5.\ Proof of Theorem~1.} \ \ It appears that Theorem~1 may be
proven by a method very similar to the one used by Krichever
\cite{Krichever1980}.

To explicitly construct meromorphic functions satisfying (\ref{e2}),
we use the basic elliptic functions $\sigma(z)$ and $\zeta(z)$ on $X =
\C/\Lambda = \C/\{\Z e_1 + \Z e_2\}$ defined by
$$
\sigma(z) =  z \prod_{m,n \in \Z^2 \setminus\{0\}} \left(1 - \frac{z}{e_{mn}}\right) \exp\left(\frac{z}{e_{mn}} +
\frac{1}{2}\frac{z^2}{e_{mn}^2}\right), \ \ e_{mn} = m e_1 + n e_2,
$$
$$
\zeta(z) = \frac{d}{dz} \log \sigma(z).
$$
It is known that
$$
\sigma(z+e_j) = -\sigma(z)e^{\eta_j(z+ e_j/2)}, \ \ \zeta(z+e_j) = \zeta(z) + \eta_j, \ \ \ j=1,2,
$$
where $\eta_j=2\zeta(e_j/2)$ satisfy $\eta_1 e_2 - \eta_2 e_1 = 2\pi
i$. Following \cite{Krichever1980}, we define
$$
\Phi(z,\alpha) = \frac{\sigma(\alpha-z)}{\sigma(\alpha)\sigma(z)}e^{\zeta(\alpha)z}, \ \ \alpha \in \C\backslash \Lambda
$$
which is $\Lambda$-periodic in $\alpha$, so that it's well defined for
$\alpha \in X\backslash\{0\}$ (this is in fact the Baker--Akhiezer
function for the one--gap Lam\'e potential).

For given multipliers $\nu_1 = e^{2\pi i \langle k, e_1 \rangle}, \nu_2 =  e^{2\pi i \langle k,e_2\rangle}$
(see (\ref{e2})) which are not of the form $e^{\beta e_1}, e^{\beta e_2}$
there are unique $\alpha \in X\backslash\{0\}$ and
$\mu\in \C$ such that any meromorphic function $\psi$ with simple
poles at $z_0 + \Lambda$ for which the periodicity conditions
(\ref{e2}) hold coincides (up to scale) with
$$
\Psi_{\mu,\alpha}(z-z_0) = e^{\mu z} \Phi(z-z_0,\alpha).
$$

Riemann--Roch implies that, for given $\nu_1$, $\nu_2$ and corresponding
$\alpha$, $\mu$, any solution to the linear problem
(\ref{e1})--(\ref{e2}) with first order poles at the $p_1$, ..., $p_n$
has to be of the form
\begin{equation}
\label{form1}
\psi = \sum_{l=1}^N a_l \Psi_{\mu,\alpha}(z-p_l).
\end{equation}
The conditions (\ref{e3}) are then fulfilled if and only if
$$
\mu a_l + \sum_{m \neq l} a_m \Phi(p_l-p_m,\alpha) = 0, \ \ \ \ l=1,\dots,N.
$$
  This
linear system has a nontrivial solution $a = (a_1,\dots,a_N) \in \C^N$
if and only if its determinant vanishes. This condition
\begin{equation}
\label{spectrum} \det A(\mu,\alpha) = 0, \ \ \
A_{lm} = \delta_{lm}\mu + (1-\delta_{lm}) \Phi(p_m-p_l,\alpha)
\end{equation}
yields an equation of the form
\begin{equation}
\label{Ncurve}
\mu^N + q_1(\alpha)\mu^{N-1} + \dots + q_N(\alpha) = 0.
\end{equation}
The spectral curve $\Gamma$ is given by the solutions to this
equation.  It is an $N$-sheeted covering of $X\backslash\{0\}$ with
the natural projection $(\mu,\alpha) \in \Gamma \to \alpha \in
X\backslash\{0\}$. The solutions $\psi$ to (\ref{e1})--(\ref{e3})
parametrized by $\Gamma$ have Floquet multipliers
\begin{equation}
  \label{eq:multiplier}
  \nu_1=\exp((\mu+\zeta(\alpha))e_1-\alpha \eta_1) \quad \textrm{and} \quad  \nu_2=\exp((\mu+\zeta(\alpha))e_2-\alpha \eta_2).
\end{equation}

Note that precisely the same condition \eqref{spectrum} (although in a
different context) appears in \cite{Krichever1980}, see (22) and the first half of (16) there assuming that $\dot{x}^k=0$ for all $k$.  In (25) of \cite{Krichever1980} it is shown that near
$\alpha=0$ the curve $\Gamma$ has one sheet at which
$\mu+\zeta(\alpha)$ has a pole at $\alpha=0$ and $N-1$ sheets at which
the Floquet multipliers \eqref{eq:multiplier} converge to multipliers
of the form $(\nu_1,\nu_2)=(e^{\beta e_1}, e^{\beta e_2})$ with
$\beta$ equal to the corresponding limit of $\mu+\zeta(\alpha)$ when
$\alpha$ goes to $0$.

Floquet functions $\psi$ for multipliers of the form
$(\nu_1,\nu_2)=(e^{\beta e_1}, e^{\beta e_2})$ cannot be obtained from
the ansatz \eqref{form1}, but from
\begin{equation}
\psi = e^{\beta z}(a_0 + \sum_{l=1}^N a_l \zeta(z-p_l))
\label{form2}
\end{equation}
with
\begin{equation}
\label{r1}
\sum_{l=1}^N a_l = 0.
\end{equation}
Condition (\ref{e3}) for this ansatz takes the form
$$
a_0 + \beta a_k + \sum_{l \neq k} a_l \zeta (p_k-p_l) = 0, \ \ \ k=1,\dots,N.
$$
After elimination of $a_0$, these conditions together with \eqref{r1}
yield $N$ linear equations for $a_1,...,a_N$. Taking the determinant
yields a polynomial equation of order $N-1$ for $\beta$. Its solutions
correspond to the $N-1$ ends of $\Gamma$ for which the Floquet
mulipliers have a limit when $\alpha$ goes to $0$.

{\sc Example.} In the simplest case $N=1$ we have
$\Gamma=X\backslash\{0\}$, the spectral curve of the one--gap
Lam\'e potential.

\medskip

{\sc Acknowledgement.} The first author (C.B.) was supported by DFG Sfb/Tr 71 ``Geometric Partial
Differential Equations'', the second author  (I.A.T.) was supported by the program of the Presidium of RAS
``Fundamental Problems of Nonlinear Dynamics in Mathematical and Physical Sciences'', in addition both authors
were supported by Hausdorff Institute of Mathematics in Bonn.

\end{document}